\documentclass[11pt,a4paper]{article}
\usepackage{amsmath, amsfonts,amsthm,amssymb,amsbsy,upref,color,graphicx,amscd,hyperref,makeidx,enumerate}
\usepackage[active]{srcltx}
\usepackage[latin1]{inputenc}

\def\Lb{\Lambda}
\def\Dr{D}
\def\g{\gamma}
\def\ra{\rightarrow}
\def\ds{\displaystyle}
\def\eps{{\varepsilon}}
\def\N{\mathbb{N}}
\def\O{\Omega}
\def\Om{\Omega}
\def\lb{\lambda}
\def\R{\mathbb{R}}
\def\A{\mathcal{A}}

\def\HH{\mathcal{H}}
\def\LL{\mathcal{L}}
\def\M{\mathcal{M}}

\newcommand{\be}{\begin{equation}}
\newcommand{\ee}{\end{equation}}
\newcommand{\bib}[4]{\bibitem{#1}{\sc#2: }{\it#3. }{#4.}}

\newcommand{\cp}{\mathop{\rm cap}\nolimits}

\numberwithin{equation}{section}
\theoremstyle{plain}
\newtheorem{teo}{Theorem}[section]
\newtheorem{lemma}[teo]{Lemma}

\newtheorem{prop}[teo]{Proposition}
\newtheorem{deff}[teo]{Definition}
\theoremstyle{remark}
\newtheorem{oss}[teo]{Remark}

\newenvironment{ack}{{\bf Acknowledgements.}}

\title{Shape Optimization Problems with Internal Constraint}

\author{Dorin Bucur, Giuseppe Buttazzo, Bozhidar Velichkov}

\begin{document}
\maketitle

\begin{abstract}
We consider shape optimization problems with internal inclusion constraints, of the form
$$\min\big\{J(\Omega)\ :\ \Dr\subset\Omega\subset\R^d,\ |\Omega|=m\big\},$$
where the set $\Dr$ is fixed, possibly unbounded, and $J$ depends on $\Omega$ via the spectrum of the Dirichlet Laplacian. We analyze the existence of a solution and its qualitative properties, and rise some open questions.

\end{abstract}



\section{Introduction}\label{sintro}

A shape optimization problem is a minimization problem of the form
\be\label{minpb}
\min\big\{J(\Omega)\ :\ \Omega\in\A\big\}
\ee
where $J$ is a suitable cost functional, possibly depending on the spectrum of an elliptic operator on $\Omega$ (in this case we speak of {\it spectral optimization problems}), and $\A$ is a  class of admissible domains. A wide literature on the subject is available, dealing with existence, regularity, necessary conditions of optimality, relaxation, explicit solutions and numerical computations of the optimal shapes. We quote for instance the books \cite{bubu05,hen06,hepi05}, where the reader may find a complete list of references on the field.

The simplest situation for the existence of a solution of problem \eqref{minpb} occurs when the class of admissible domains $\A$ satisfies an external inclusion constraint, i.e. consists on quasi-open sets which are supposed {\it a priori} contained in a given {\it bounded} open set $D$ of the Euclidean space $\R^d$, 
$$\A=\{\Omega\ :\ \Omega\subset\Dr,\ \Omega\mbox{ quasi-open}\}.$$
In this case a general existence result, due to Buttazzo and Dal Maso (see \cite{budm93}), states that problem \eqref{minpb}, with the additional constraint $|\Omega|\le m$ on the Lebesgue measure of the competing domains, admits a solution provided the cost functional $J$ satisfies some mild conditions:

\begin{itemize}
\item[(i)]$J$ is lower semicontinuous for the $\gamma$-convergence, suitably defined;
\item[(ii)]$J$ is monotone decreasing for the set inclusion.
\end{itemize}

\noindent Some interesting cases fall into the framework above, as for instance the ones below.

\medskip{\it Spectral optimization.} For every admissible domain $\Omega$ consider the Dirichlet Laplacian $-\Delta$ which, under mild conditions on $\Omega$, admits a compact resolvent and so a discrete spectrum $\lambda(\Omega)$. The cost is in this case of the form
$$J(\Omega)=\Phi\big(\lambda(\Omega)\big)$$
for a suitable function $\Phi$. For instance, taking $\Phi(\lambda)=\lambda_k$ we may consider the optimization problem for the $k$-th eigenvalue of $-\Delta$:
$$\min\big\{\lambda_k(\Omega)\ :\ \Omega\in\A\big\}.$$

\medskip{\it Integral functionals.} Given a right-hand side $f$ we consider the PDE
$$-\Delta u=f\hbox{ in }\Omega,\qquad u\in H^1_0(\Omega)$$
which provides, for every admissible domain $\Omega$, a unique solution $u_\Omega$ that we assume extended by zero outside of $\Omega$. The cost is in this case of the form
$$J(\Omega)=\int_{\R^d}j\big(x,u_\Omega(x)\big)\,dx$$
where $j$ is a given integrand.

\medskip When the surrounding box $D$ is unbounded the existence result above is no longer true, as some simple examples show. In the case $D=\R^d$ a quite different approach to the proof of the existence of optimal domains has been considered by Bucur in \cite{buc00}, using a refined argument related to the Lions concentration-compactness principle (see \cite{pll84}).

In this paper we consider problem \eqref{minpb} where the admissible class $\A$ is defined through an internal constraint:
\be\label{admco}
\A=\big\{\Omega\ :\ \Dr\subset\Omega\subset\R^d,\ \Omega\mbox{ quasi-open}, \ |\Omega|\le m\big\}
\ee
where $\Dr$ is a fixed quasi-open set of finite measure, possibly unbounded. 

In spite of its simplicity, even for cost functionals like $J(\Omega)=\lambda_1(\Omega)$, the existence proof is rather involved, and several interesting questions arise. For this functional, together with the existence of a solution, we prove some global properties for the optimal set: it has to lie in finite distance to $\Dr$ (in particular the optimal set is bounded, provided $\Dr$ is bounded), it has finite perimeter outside $\overline\Dr$, it is an open set as soon as its measure is strictly greater than the measure of the (quasi-connected) $\Dr$. Local regularity properties, outside $\overline \Dr$ are not discussed here, being similar to the bounding box situation, and we refer the reader for instance to \cite{brla}. We discuss as well the existence question for $J(\Omega)=\lambda_k(\Omega)$, and refer the reader to \cite{bulbk} for the analysis of these functionals in the absence of any inclusion constraint in $\R^d$.

\section{Notations and preliminaries}\label{sprel}

We introduce here the main tools we use in the following; further details can be found for instance in \cite{bubu05,but10}.

In the sequel, we will work in the Euclidean space $\R^d$ with $d\ge2$. Given a subset $E\subset\R^d$ we define the capacity of $E$ by
$$\cp(E)=\inf\Big\{\|u\|^2_{H^1}\ :\ u\in{\cal U}_E\Big\}\,,$$
where ${\cal U}_E$ is the set of all functions $u$ of the Sobolev space $H^1(\R^d)$ such that $u\ge1$ almost everywhere in a neighborhood of $E$. If a property $P(x)$ holds for all $x\in E$ except for the elements of a set $Z\subset E$ with $\cp(Z)=0$, we say that $P(x)$ holds {\it quasi-everywhere} (shortly {\it q.e.}) on $E$, whereas the expression {\it almost everywhere} (shortly {\it a.e.}) refers, as usual, to the Lebesgue measure, that we often denote by $|\cdot|$.

A subset $\Omega$ of $\R^d$ is said to be {\it quasi-open} if for every $\eps>0$ there exists an open subset $\Omega_\eps$ of $\R^d$, with $\Omega\subset\Omega_\eps$, such that $\cp(\Omega_\eps\setminus\Omega)<\eps$. Similarly, a function $f:\R^d\to\R$ is said to be {\it quasi-continuous} (resp. {\it quasi-lower semicontinuous}) if there exists a decreasing sequence of open sets $(\omega_n)_{n>0}$ such that $\lim_{n\to\infty}\cp{\omega_n}=0$ and the restriction $f_n$ of $f$ to the set $\omega_n^c$ is continuous (resp. lower semicontinuous). It is well known (see for instance \cite{zie89}) that every function $u\in H^1(\R^d)$ has a quasi-continuous representative $\tilde u$, which is uniquely defined up to a set of capacity zero, and given by
$$\tilde u(x)=\lim_{\eps\to0}\frac{1}{|B_\eps(x)|}\int_{B_\eps(x)}u(y)\,dy\,,$$
where $B_\eps(x)$ denotes the ball of radius $\eps$ centered at $x$. We often identify the function $u$ with its quasi-continuous representative $\tilde u$; in this way, we have that quasi-open sets can be characterized as the sets of strict positivity of functions in $H^1(\R^d)$ and that the capacity can be equivalently defined by
$$\cp(E)=\min\Big\{\|u\|^2_{H^1}\ :\ u\in H^1(\R^d),\ u\ge1\hbox{ q.e. on }E\Big\}.$$

For every quasi-open set $\Omega\subset\R^d$ we denote by $H^1_0(\Omega)$ the space of all functions $u\in H^1(\R^d)$ such that $u=0$ q.e. on $\R^d\setminus\Omega$, with the Hilbert space structure inherited from $H^1(\R^d)$,
$$\langle u,v\rangle_{H^1_0(\Omega)}=\langle u,v\rangle_{H^1(\R^d)}.$$
The usual properties of Sobolev functions on open sets extend to quasi-open sets.

Let $\Omega$ be a quasi-open set of finite measure. By $R_\Omega$ we denote the resolvent operator of the Laplace equation with Dirichlet boundary condition,
$$R_\Omega:L^2(\R^d)\to L^2(\R^d),$$
where $R_\Omega(f)$ is the weak solution of the equation
$$\begin{cases}
-\Delta u=f\in L^2(\R^d),\\
u\in H^1_0(\Omega).
\end{cases}$$

We denote by $\M_0$ the set of capacitary measures on $\R^d$, that is the set of all Borel measures, possibly taking the value $+\infty$, vanishing on all sets of zero capacity. Observe that for each Borel set $S$ the measure
$$\infty_S(B)=\begin{cases}
0&\hbox{if }\cp(B\cap S)=0\\
+\infty&\hbox{otherwise}
\end{cases}$$
is a capacitary measure.

For each capacitary measure $\mu$, we define the linear vector space
$$H^1(\R^d)\cap L^2(\R^d,\mu)=\Big\{u\in H^1(\R^d)\ :\ \int_{\R^d}|u|^2\,d\mu<\infty\Big\}.$$
Taking $\mu=\infty_S$ with $S=\Omega^c$ gives $H^1(\R^d)\cap L^2(\R^d,\mu)=H^1_0(\Omega)$. In \cite{budm91} it was shown that the above space, endowed with the scalar product
$$\langle u,v\rangle=\int_{\R^d}\nabla u\nabla v\,dx+\int_{\R^d}uv\,dx+\int_{\R^d}uv\,d\mu,$$
is a Hilbert space. Moreover, the space $H^1(\R^d)\cap L^2(\R^d,\mu)$ is separable when seen as a subset of the separable metric space $H^1(\R^d)$. If $\{u_n\}_{n\ge0}\subset H^1(\R^d)\cap L^2(\R^d,\mu)$ is a dense countable subset, then we define the regular set of the capacitary measure $\mu\in\M_0$ as
$$\Omega_\mu=\bigcup_{n\ge0}\{u_n\ne0\}.$$
If the set $\Omega_{\mu}$ has finite Lebesgue measure, then
$$\|u\|^2=\int_{\R^d}|\nabla u|^2\,dx+\int_{\R^d}|u|^2\,d\mu,$$
is an equivalent norm on $H^1(\R^d)\cap L^2(\R^d,\mu)$.
We define the resolvent $R_\mu$ as the map
$$R_\mu:L^2(\R^d)\to L^2(\R^d),$$
which associates to each function $f\in L^2(\R^d)$ the solution $u$ of the relaxed problem formally written as
$$-\Delta u+\mu u=f,\qquad u\in H^1(\R^d)\cap L^2(\R^d,\mu),$$
which has to be rigorously defined in the weak form
$$\begin{cases}
\ds\int_{\R^d}\nabla u\nabla\varphi\,dx+\int_{\R^d}u\varphi\,d\mu
=\int_{\O}f\varphi\,dx\qquad\forall\varphi\in H^1(\R^d)\cap L^2(\R^d,\mu),\\
u\in H^1(\R^d)\cap L^2(\R^d,\mu). 
\end{cases}$$

If $\mu$ is a capacitary measure with regular set of finite Lebesgue measure, then the constant function $1$ is in the dual space of $H^1(\R^d)\cap L^2(\R^d,\mu)$, so we can define $w_\mu:=R_\mu(1)$ and we have $\Om_\mu=\{w_\mu>0\}$ up to zero capacity sets.

We consider the following relation of equivalence on $\M_0$:
$$\mu_1\sim\mu_2\ \iff\ \mu_{1}(\Omega)=\mu_2(\Omega),\qquad\forall\Omega\mbox{ quasi-open},$$
and we can make the quotient set $\M_0/\sim$ a metric space (we still denote this quotient by $\M_0$ and call its elements capacitary measures), by introducing the convergence below.

\begin{deff}
To each bounded open set $\Omega$ and each capacitary measure $\mu$ we associate the functional
$$F_\mu(u,\Omega)=\int_{\R^d}|\nabla u|^2\,dx+\int_{\R^d}u^2\,d\mu+\chi_{H^1_0(\Omega)}(u),$$
defined on the metric space $L^2(\R^d)$, where
$$\chi_{H^1_0(\Omega)}(u)=
\begin{cases}
0&\hbox{if }u\in H^1_0(\Omega),\\ 
+\infty&\hbox{otherwise.}
\end{cases}$$
We say that a sequence $(\mu_n)_{n\ge0}\subset\M_0$ $\gamma_{loc}$-converges (locally $\gamma$-converges) to $\mu\in\M_0$ and we write
$$\mu_n\stackrel{\gamma_{loc}}{\to}\mu,$$
if for each bounded open set $\Omega\subset\R^d$, 
$$F_{\mu_n}(\cdot,\Omega)\stackrel{\Gamma}{\to}F_\mu(\cdot,\Omega),$$
where the above expression denotes the usual $\Gamma$ convergence of functionals on the metric space $L^2(\R^d)$, that is
\begin{enumerate}
\item[(i)]for every $u\in L^2(\R^d)$ and every sequence $(u_n)_{n\ge0}\subset L^2(\R^d)$ converging to $u$ in the norm of $L^2(\R^d)$ we have
$$F_\mu(u,\Omega)\le\liminf_{n\to\infty}F_{\mu_n}(u_n,\Omega),$$
\item[(ii)]for every $u\in L^2(\R^d)$ there exists a sequence $(u_n)_{n\ge0}\subset L^2(\R^d)$ converging to $u$ in the norm of $L^2(\R^d)$ such that
$$F_\mu(u,\Omega)=\lim_{n\to\infty}F_{\mu_n}(u_n,\Omega).$$
\end{enumerate}
\end{deff}

\begin{oss}
In \cite[Definition 2.7]{badm87} the $\g_{loc}$-convergence introduced above was called $\gamma$-convergence (see also \cite{dmmo87}). Here, we chose to denote by $\gamma$-convergence a stronger convergence, as follows.
\end{oss}

\begin{deff}
Let $({\mu_n})_{n\ge0}$ and $\mu$ be capacitary measures such that their regular sets have uniformly bounded Lebesgue measures. We say that $\mu_n$ $\gamma$-converges to $\mu$, if $(w_{\mu_n})_{n\ge0}$ converges in $L^{2}(\R^d)$ to $w_{\mu}$.
\end{deff}
\begin{oss}\label{remlam}
With the definition above, we have the equivalence
$$\mu_n\stackrel{\gamma}{\longrightarrow}\mu\quad\iff\quad R_{\mu_n}\stackrel{\LL(L^2(\R^d))}{\longrightarrow}R_\mu.$$
Indeed, for the "$\Rightarrow$" implication, we refer to \cite[Proposition 3.3]{buc00}. For the converse implication, the proof is immediate. On the one hand, we have
$$R_{\mu_n}(1_{\Om_{\mu_n}})- R_{\mu}(1_{\Om_{\mu_n}}) \rightarrow 0\;\;\mbox { in}\;\;L^2(\R^d),$$
and on the other hand
$$R_{\mu_n}(1_{\Om_{\mu}})-R_{\mu}(1_{\Om_{\mu}})\rightarrow 0\;\;\mbox { in}\;\;L^2(\R^d).$$
Making the difference we get that
$$\|R_{\mu_n}(1_{\Om_{\mu_n}})-R_{\mu_n}(1_{\Om_{\mu}})+R_{\mu}(1_{\Om_{\mu}})-R_{\mu}(1_{\Om_{\mu_n}})\|_{L^2(\R^d)} \rightarrow 0$$
 and using the maximum principle we conclude with
 $$\|R_{\mu_n}(1)-R_{\mu}(1)\|_{L^2(\R^d)} \rightarrow 0.$$
\end{oss}

\begin{oss}
In the case $\mu=\infty_{A^c}$ with $A$ quasi-open, we have 
$$F_A(u,\Omega)=F_{\infty_{A^c}}(u,\Omega)=\int_{\R^d}|\nabla u|^2\,dx+\chi_{H^1_0(\Omega\cap A)}(u).$$
\end{oss}

The following compactness theorem was proved in \cite{dmmo87}.

\begin{teo}\label{comgam}
For every sequence $(\mu_n)_{n\ge0}\subset\M_0$ there exists a subsequence $(\mu_{n_k})_{k\ge0}$ which $\gamma_{loc}$-converges to a measure $\mu$ of the class $\M_0$.
\end{teo}

In \cite{buc00} it was shown that if a sequence of quasi-open sets $(\Omega_n)_{n\ge0}$ of uniformly bounded measure $\gamma_{loc}$-converges to a capacitary measure $\mu$, then the sequence of functionals $F_{\Omega_n}=F_{\infty_{\Omega_n^c}}$ $\Gamma$-converges in $L^2(\R^d)$ to the functional $F_\mu$, where
$$F_{\Omega_n}(u)=\int_{\R^d}|\nabla u|^2\,dx+\chi_{H^1_0(\Omega_n)}(u),$$ 
$$F_\mu(u)=\int_{\R^d}|\nabla u|^2\,dx+\int_{\R^d}u^2\,d\mu.$$
Furthermore, in the same work, Theorem 5.4 states the following.

\begin{teo}
If $(\Omega_n)_{n\ge0}$ is a sequence of quasi-open sets of uniformly bounded measure which $\gamma_{loc}$-converges to a capacitary measure $\mu$ and is such that the sequence of solutions $(w_{\Omega_n})_{n\ge0}$ converges in $L^2(\R^d)$ to some function $w$, then $\{w>0\}$ is the regular set of $\mu$ and $w=w_\mu$.
\end{teo}

\begin{oss}
The result above is actually valid (with practically the same proof) for sequences of measures whose regular sets have uniformly bounded Lebesgue measure. From this observation and the fact that the $\gamma_{loc}$-convergence is metrizable (see \cite{dmmo87}, Proposition 4.9) with metric $d_{\gamma_{loc}}$, we obtain that for each $t>0$ the space
$$\M_0^t=\{\mu\in\M_0\ :\ |\Omega_\mu|\le t\},$$
where $\Omega_\mu$ denotes the regular set of $\mu$, is a complete metric space when endowed with the metric
$$d(\mu_1,\mu_2)=d_{\gamma_{loc}}(\mu_1,\mu_2)+\|w_{\mu_1}-w_{\mu_2}\|_{L^2(\R^d)}.$$
Moreover, suppose that $(\mu_n)_{n\ge0}$ is a sequence such that $(w_{\mu_n})_{n\ge0}$ converges in $L^2(\R^d)$. Then each subsequence of $(\mu_n)_{n\ge0}$ has a convergent subsequence in the metric $d$ and the limit is uniquely determined by the limit $w=L^2$-$\lim_{n\to\infty}w_{\mu_n}$. Since we are in a metric space, we have that $(\mu_n)_{n\ge0}$ converges to some capacitary measure $\mu\in\M_0^t$ in the metric $d$. Then, the space $\M_0^t$ endowed with the metric $d_{\gamma}$, where 
$$d_{\gamma}(\mu_1,\mu_2)=\|w_{\mu_1}-w_{\mu_2}\|_{L^2(\R^d)},$$
is a complete metric space.
\end{oss}

\begin{lemma}\label{lscig}
Consider a sequence $(\O_n)_{n\ge0}$ of quasi-open sets of uniformly bounded measure such that $\O_n$ $\gamma$-converges to the capacitary measure $\mu$ with regular set $\O_\mu$. Then, for every $k\ge 1$
$$\lambda_k(\O_\mu)\le\lambda_k(\mu)=\lim_{n\to\infty}\lambda_k(\O_n).$$
\end{lemma}

\begin{proof}
By Remark \ref{remlam}, $R_{\O_n}\to R_\mu$ in the operator norm of $\LL(L^2(\R^d))$, and so we have
$$\lb_k(\O_n) \ra \lb_k(\mu).$$
The inequality 
$$\lambda_k(\O_\mu)\le\lambda_k(\mu),$$
now follows as a consequence of the inequality of the measures $\infty_{\O_\mu^c}(B)\le\mu(B)$, for each quasi-open set $B$, in the $\min$-$\max$ definition of the eigenvalues.
\end{proof}

\section{Existence of an optimal set}\label{stheo}

For $\A$ defined in \eqref{admco}, we study the existence of a solution for the problem
\be\label{lbk}
\min\{\lambda_k(\O)\ :\ \O\in\A\}.
\ee
We notice that if a solution $\O$ of problem \eqref{lbk} exists, then necessarily the measure of $\O$ is precisely equal to $m$. Indeed, assume by contradiction that $\O$ is an optimal set with measure strictly less than $m$. There exists an open set $U$ with measure still less than $m$ which contains $\O$, so it is an optimal set too. Since $U$ is open, we can add small balls on each connected component, so that the global measure still remains less than $m$, but the $k$-th eigenvalue strictly diminishes, which contradicts the optimality of $\Om$.

The fundamental tool, allowing to understand the behaviour of a minimizing sequence in $\R^d$, is the following concentration-compactness result (see \cite[Theorem 2.2]{buc00}) for the resolvent operators.

\begin{teo}\label{thcc}
Let $(\O_n)_{n\ge0}$ be a sequence of quasi-open sets of uniformly bounded measure. Then, there exists a subsequence, still denoted by $(\Omega_n)_{n\ge0}$, such that one of the following situations occurs.
\begin{enumerate}[(i)]
\item{\textbf{Compactness.}} There exists a sequence of vectors $(y_n)_{n\ge0}\subset\R^d$ and a capacitary measure $\mu$, such that $y_n+\O_n$ $\gamma$-converges to the measure $\mu$, and so $R_{y_n+\O_n}$ converges in the uniform operator topology of $L^2(\R^d)$ to $R_\mu$.
\item{\textbf{Dichotomy.}} There exists a sequence of subsets $\tilde{\O}_n\subseteq \O_n$, such that:
\begin{itemize}
\item $\|R_{\O_n}-R_{\tilde{\O}_n}\|_{\LL(L^2(\R^d),L^2(\R^d))}\to0$;
\item $\tilde{\O}_n$ is a union of two disjoint quasi-open sets $\tilde{\O}_n=\O_n^+\cup \O_n^-$;
\item $d(\O_n^+,\O_n^-)\to\infty$;
\item $\liminf_{n\to\infty}|\O_n^\pm|>0$.
\end{itemize} 
\end{enumerate}
\end{teo}

\begin{oss}\label{remtrans}
Assume that $\Dr$ is a non-empty quasi-open set and $(\O_n)_{n\ge1}$ is a sequence of quasi-open sets such that $\Dr \subset \Om_n$, $|\Om_n|\le m$. Then, if the compactness situation holds in Theorem \ref{thcc}, then one can take $y_n=0$, i.e. no translation is necessary. 

Suppose first that $y_n$ is divergent and notice that the solution $w_{\Dr+y_n}$ is just $w_{\Dr}$ translated to the left by $y_n$.
By the maximum principle, we have that $w_{\O_n+y_n}\ge w_{\Dr+y_n}$ and so
$$\int w_{\Dr+y_n}w_{\O_n+y_n}\,dx\ge\int w_{\Dr}^2\,dx>0.$$
Since $y_n\to\infty$, we have that $w_{\Dr+y_n}\rightharpoonup0$ weakly in $L^2$. By the strong convergence of $w_{\O_n+y_n}$ we have
$$\int w_{\Dr+y_n}w_{\O_n+y_n}\,dx\to0,$$
which is a contradiction and so we have that $y_n$ is bounded.\\
Choose a convergent subsequence $y_{n_k}\to y$ and set $w=L^2$-$\lim_{\tau \to\infty}w_{\O_{n_k}+y_{n_k}}$. We have
$$\begin{array}{ll}
\|w_{\O_{n_k}}-w(\cdot-y)\|_{L^2(\R^d)}&\le\|w_{\O_{n_k}}-w(\cdot-y_{n_k})\|_{L^2(\R^d)}+\|w(\cdot-y_{n_k})-w(\cdot-y)\|_{L^2(\R^d)}\\
&\le\|w_{\O_{n_k}+y_{n_k}}-w\|_{L^2(\R^d)}+\|w(\cdot-y_{n_k})-w(\cdot-y)\|_{L^2(\R^d)},
\end{array}$$
and both last terms converge to zero as $k\to\infty$. By extracting another subsequence we obtain a subsequence converging in $L^2(\R^d)$ to the function $w(\cdot-y)$ and $\gamma_{loc}$-converging to a capacitary measure $\mu$. By \cite[Theorem 5.4]{buc00}, we obtain $w=w_\mu$.
\end{oss}

\begin{teo}\label{mth1}
Let $\Dr$ be a quasi-open set of finite measure. Then, the problem
\be\label{lbk1}
\min\{\lambda_1(\O)\ :\ \O\mbox{ quasi-open},\ \Dr\subset\O,\ |\O|\le m\}
\ee
has at least one solution.
\end{teo}

\begin{proof}
We consider a minimizing sequence $(\O_n)_{n\ge1}$ with the property that $\liminf_{n\to\infty}|\Om_n|$ is minimal. Clearly, this value can not be equal to zero. According to Theorem \ref{thcc} and Remark \ref{remtrans}, if we are in the compactness situation, for a subsequence (still denoted with the same indices) there exists a measure $\mu$ such that $\Om_n$ $\g$-converges to $\mu$ and 
$R_{\Om_n}$ converges to $R_\mu$ in the uniform operator topology of $L^2(\R^d)$. As a consequence, by Lemma \ref{lscig} we get that the regular set $\Omega_\mu$ of $\mu$ is a solution.

If we are in the dichotomy situation, we get a contradiction. On the one hand since $\Om_n^+$ and $\O_n^-$ are at positive distance, we may assume that $\lb_1(\O_n^+\cup\O_n^-)=\lb_1(\O_n^+)$. Then, the sequence $\O_n^+\cup\Dr$ is also minimizing since $|\lb_1(\O_n^+)-\lb_1(\Om_n)|\to0$ (see \cite[Proposition 3.7]{buc00}), but either
$$\liminf_{n\to\infty}|\O_n^+\cup\Dr|<\liminf_{n\to\infty}|\O_n|$$
or $|\O_n^-\setminus\Dr|\to0$. The first assertion is in contradiction with our assumption on the choice of a least measure minimizing sequence. The second assertion is also impossible, since it implies that $d(\O_n^+,\{0\})\to+\infty$, otherwise the measure of $\Dr$ would be infinite. Consequently, since the measure of $\Dr$ is finite, we get that $|\O_n^+\cap\Dr|\to0$ and consider the ball $B$ of measure equal to $\limsup |\O_n^+|$. Therefore, $B \cup \Dr$ is a solution for every position of the ball B. In particular, this leads to a contradiction if the ball intersects, but not cover, a quasi-connected component of $\Dr$.
\end{proof}

\begin{oss}
Let us notice that from every minimizing sequence we can extract a $\g$-convergent subsequence. The basic observation is that any minimizing sequence for which $\liminf_{n\to\infty}|\Om_n|$ is minimal leads to an optimal set, which necessarily has the measure equal to $m$. Since the measure is lower semicontinuous for the $\g$-convergence, this means that any minimizing sequence should satisfy $\lim_{n\to\infty}|\Om_n|=m$ excluding the dichotomy in the proof above. 
\end{oss}

In the sequel we show a result which gives a rather explicit behavior of a minimizing sequence for the problem
\be\label{lbkk}
\min\big\{\lambda_k(\O)\ :\ \O\mbox{ quasi-open},\ \Dr\subset\O,\ |\O|\le m\big\}.
\ee
For every $m>0$ we introduce the value
$$\lb_k^*(m)=\inf\big\{\lambda_k(\O)\ :\ \O\mbox{ quasi-open},\ |\O|\le m\big\}.$$
Following \cite{bulbk}, there exists a bounded quasi-open set $\O$, with measure equal to $m$ such that $\lb_k(\O)=\lb_k^*(m)$.

\begin{teo}\label{mth2}
For $k\in\N$, $k \ge 2$, let 
$$\alpha_k=\inf\{\lambda_k(\O)\ :\ \O\mbox{ quasi-open},\ \Dr\subset\O,\ |\O|\le m\}.$$
One of the following assertions holds:
\begin{itemize}
\item[(i)]problem \eqref{lbkk} has a solution;
\item[(ii)]there exists $l\in\{1,\dots,k-1\}$ and an admissible quasi-open set $\O$ such that $\alpha_k=\lb_{k-l}(\O)=\lb_l^*(m-|\O|)$;
\item[(iii)]there exists $l\in \{1,\dots,k-1\}$ such that $\alpha_k=\lb_l^*(m-|\O|)> \lb_{k-l}(\Dr)$.
\end{itemize}
\end{teo}

\begin{proof}
Let us consider a minimizing sequence $(\O_n)_{n\ge1}$ with the property that $\liminf_{n\to\infty}|\Om_n|$ is minimal. If compactness occurs in Theorem \ref{thcc}, then the existence of a solution follows as in Theorem \ref{mth1}.

If dichotomy occurs, as in Theorem \ref{mth1} we may assume that
$$|\Om_n^+|\to\alpha^+,\qquad|\Om_n^-|\to\alpha^-,\qquad|\Om_n^+\cap \Dr|\to0.$$
Then, up to a subsequence there exists $l\in\{1,\dots,k-1\}$ such that one of the two possibilities below holds:
\begin{itemize}
\item [(A)]\quad$|\lb_k(\Om_n)-\lb_{k-l}(\Om_n^-)|\to0$\quad and\quad$\lb_l(\Om_n^+)\le\lb_{k-l}(\Om_n^-)\le\lb_{l+1}(\Om_n^+)$;
\item [(B)]\quad$|\lb_k(\Om_n)-\lb_{l}(\Om_n^+)|\to0$\quad and\quad$\lb_{k-l}(\Om_n^-)\le\lb_l(\Om_n^+)\le\lb_{k-l+1}(\Om_n^-)$.
\end{itemize}
We may take the maximal $l$ with such a property. We use now an induction argument as follows. For $k=1$ as proved in Theorem \ref{mth1}, dichotomy does not occur, so the compactness gives (i). Assume that for $1,\dots,k-1$ Theorem \ref{mth2} is true. We prove it for $k$. If compactness occurs, then (i) holds. If dichotomy occurs and we are in situation (A) we get that $(\O_n^-\cup\Dr)_n$ is minimizing for the $k-l$ eigenvalue with the inclusion constraint and the corresponding measure $m-\alpha^+\ge\liminf |\O_n^-\cup\Dr|$. Since $l$ is maximal with this property, for the sequence $(\O_n^-\cup\Dr)_n$ dichotomy cannot occur again, so finally (ii) holds.

If (B) occurs, then $|\Om_n^-\setminus\Dr|\to0$ and we are in situation (iii).
\end{proof}

\begin{oss}
Theorem \ref{mth2} gives a complete description on the behaviour of a minimizing sequence for $\lb_k$, $k \ge 2$. Assertion i) implies the existence of a solution. As well, if $\Dr$ has some suitable geometric properties, both alternatives (ii) and (iii) lead to the existence of a solution. Typically, if for every $R>0$, there exists $x\in \R^d$ such that $B_R(x)\cap D=\emptyset$, then existence of a solution occurs. The main ingredients for obtaining such a result are the existence of a minimizer for the $k$-th eigenvalue in absence of any inclusion constraint (see \cite{bulbk}) and the boundedness result obtained in Proposition \ref{boundedness} of the next section and in \cite{bulbk}. We do not know whether existence holds without this assumption on $\Dr$, but we were not able to find a counterexample (see the last section). 
\end{oss}

\section{Qualitative properties of the optimal sets}

A natural question that arises in the shape optimization problems with constraints like \eqref{lbk1} is to understand the influence of the inclusion domain $\Dr$ on the optimal sets: does boundedness and/or convexity of $\Dr$ imply the same properties on the optimal set? As we shall see, the answer is positive for the boundedness constraint, but negative for the convexity constraint.

\subsection{Regularity of the optimal set}

In this section we deal with the penalized version of problem \eqref{lbk1} 
\be\label{sop}
\min\big\{\lambda_1(\Omega)+\Lb|\Omega|:\Omega\subset\R^d,\ \Omega\hbox{ quasi-open, }\Dr\subset\Omega\big\},
\ee
for some $\Lb>0$. For the local equivalence of the two problems we refer the reader to \cite{brla}. As well, we refer the reader to \cite{brla} for a complete analysis of a similar problem, in which the internal constraint $\Dr\subset\Omega$ is replaced by an external constraint $\Om \subset \Dr$, with a bounded open set $\Dr$.

In this section we will analyze the internal constrain problem, and prove that the optimal set of \eqref{sop} is open even if $\Dr$ is only quasi-open, provided that $\Dr$ is quasi-connected and the optimal set has a measure strictly greater than $ |\Dr|$. 
For simplicity, we say that $\Dr$ is quasi-connected if for every couple of non-empty quasi-open sets $A_1$ and $A_2$ having intersection of positive capacity with $\Dr$ and such that $\Dr \subset A_1\cup A_2$, we get $\cp (A_1\cap A_2)>0$.
The quasi-connectedness has a topological counterpart. Indeed, a quasi-open, quasi-connected set $A$ has a fine interior (which differs from $A$ by a set of zero capacity) which is finely connected (the fine topology being the coarsest topology making all superharmonic functions continuous). A nonnegative superharmonic function in $H_0^1(A)$ with $A$ finely connected, is either equal to $0$ or is strictly positive (see \cite{bbf99,fug72,lat98}).

In the following, without loss of generality we assume that $\Lb=1$.

\begin{oss}
The existence of a solution to \eqref{sop} follows by the same argument we used in the proof of Theorem \ref{mth1} and so we omit the proof.
\end{oss}

Let $\Dr$ be a quasi-open, quasi-connected set of finite measure. Let $\Omega$ be a solution of problem \eqref{sop}, let $\lambda:=\lambda_1(\Omega)$, and let $u:=u_{\Omega}$ be the first normalized eigenfunction:
$$\begin{cases}
-\Delta u=\lambda u,\\
u\in H^1_0(\Omega),\quad\|u\|_{L^2}=1.
\end{cases}$$
As $\Dr$ is quasi-connected, if $\Omega$ is optimal, then $u$ is a solution of the optimization problem
\be\label{op}
\min\Big\{\frac{\int|\nabla v|^2\,dx}{\int v^2\,dx}+|\{v>0\}|\ :\ v\in H^1(\R^d),\ \Dr\subset\{v>0\}\Big\}.
\ee

The following Lemma is similar to \cite[Lemma 3.2]{altcaf} and \cite[Lemma 3.1]{brla}.
 
\begin{lemma}\label{boundbelow}
Let $u$ be a solution of the problem \eqref{op}. Then there is a constant $C$ depending only on the dimension $d$ such that for each $r>0$, the following implication holds:
\be
Cr\le\frac{1}{|\partial B_r|}\int_{\partial B_r}u\,d\HH^{d-1}\quad\Longrightarrow\quad B_r\subset\{u>0\}.
\ee
\end{lemma}
The next proposition follows the approach first introduced in  \cite{altcaf}; nevertheless, we give the proof below to stress the fact that the quasi-open internal constraint does not change the argument too much.

\begin{prop}\label{open}
Let $\Dr$ be a quasi-open, quasi-connected set of finite measure. Every solution $\Omega$ of problem \eqref{sop} is an open set up to a set of capacity $0$.
\end{prop}

\begin{proof}
Let $u$ be a solution of  \eqref{op}. We  prove that if $u(x)>0$, then $u$ is positive in a small ball centered at $x$. Without loss of generality, we can suppose that $x=0$ and that $0$ is a regular point  in the sense that 
$$u(0)=\lim _{r\ra 0} \frac{1}{|B_r(0)|}\int_{B_r(0)}u(y) dy.$$
 Denote by $\varphi_r$ the solution of 
\be\label{solwith1}
\begin{cases}
-\Delta\varphi_r=1,\\
\varphi_r\in H^1_0(B_r),
\end{cases}
\ee
where $B_r$ denotes the ball centered in $0$ of radius $r$. An explicit computation gives
$$\varphi_r(y)=\frac{r^2-|y|^2}{2d}\;.$$
Since
$0\le\Delta u+\lambda u$ in the distributional sense on $\R^d$, we have 
$$\Delta(u-\|u\|_{\infty}\lambda\varphi_r)\ge-\lambda u+\lambda\|u\|_{\infty}\ge0$$
on each ball $B_r$, so  the function 
$u-\|u\|_{\infty}\lambda\varphi_r$ is subharmonic on $B_r$.
By the Poisson's formula, we have
$$u(0)-\|u\|_{\infty}\lambda\varphi_r(0)\le C(d)\frac{1}{|\partial B_r|}\int_{\partial B_r}u(y)\,d\HH^{d-1}(y),$$
$$u(0)-\|u\|_{\infty}\lambda C_1 r^2\le C(d)\frac{1}{|\partial B_r|}\int_{\partial B_r}u(y)\,d\HH^{d-1}(y).$$
Suppose that $u(0)>0$. Then, choosing $r$ small enough, we have
$$u(0)\le  \frac{2 C(d)}{|\partial B_r|}\int_{\partial B_r}u(y)\,d\HH^{d-1}(y).$$
Now choose $C$ as in Lemma \ref{boundbelow} and $r$ such that $2 rCC(d)\le u(0)$. Then
$$Cr\le\frac{1}{|\partial B_r|}\int_{\partial B_r}u(y)\,d\HH^{d-1}(y),$$
and so $u>0$ on $B_r$.
\end{proof}

\begin{oss}
Alternatively, one can formulate the proposition above, requiring that the inclusion $\Dr\subset\Om$ holds quasi-everywhere, and in this case the optimal sets $\{u>0\}$ in \eqref{op} are open.
\end{oss}

\begin{oss}
If $\Dr$ is a quasi-open set such that there does not exist an open set containing $\Dr$ and having the same Lebesgue measure, then the Proposition above asserts that the measure of any optimal set is strictly greater than the measure of $\Dr$. 

In general, this is not the case if $\Dr$ is an open set. Indeed, following a simple computation one can consider $\Dr$ to be a ball $B$ and take a constant $\Lb$ large enough, so that the optimal set is $B$ itself. More generally, if the partial metric derivative of the first eigenvalue on $\Dr$ is finite, i.e.
$$\lambda'_1(\Dr):=\limsup_{|\tilde\Dr\setminus\Dr|\to0,\ \tilde\Dr\supset\Dr}\frac{\lb_1(\Dr)-\lb_1(\tilde\Dr)}{|\tilde\Dr\setminus\Dr|}<+\infty,$$
then for every $\Lb>\lb'_1(\Dr)$ there exists $\Lb'>\Lb$ such that the optimal solution of \eqref{sop} with $\Lb'$ is $\Dr$. Indeed, by contradiction for every $\Lb>\lb'_1(\Dr)$ there exists $\eps>0$ such that for every $\Omega\supset\Dr$ such that $|\O|\le|\Dr|+\eps$ we have
$$\lb_1(\Dr)+\Lb|\Dr|\le\lb_1(\O)+\Lb|\O|.$$
Then, replacing $\Lb$ with $\Lb'>\Lb$ such that $\Lb'\ge\frac{\lb_1(\Dr)}{\eps}$, we get that $\Dr$ is a global minimizer.
\end{oss}

\subsection{Bounded constraint implies bounded minimizers}

We give the following technical result for which we refer to \cite[Lemma 3.4]{altcaf} and to \cite[Lemma 3.1]{brla}.

\begin{lemma}\label{sopra}
For every solution $u$ of the optimization problem \eqref{op}, there exists a constant $C_0$ and $r_0$ such that for each $x\in\R^d$ such that $d(x,\Dr)>r_0$ and for each $r<r_0$ the following implication holds:
\be\label{boundabove}
\left(\frac{1}{|\partial B_r(x)|}\int_{\partial B_r}u\,d\HH^{d-1}\le C_0r\right)\Rightarrow\left(u=0\ on\ B_{\frac{ r}{2}}\right).
\ee
The constants $C_0$ and $r_0$ above depend only on the dimension $d$ of the space and on $\lb_1(\Dr)$ respectively.
\end{lemma}

\begin{prop}\label{boundedness}
Suppose that $\Dr$ is a quasi-open set of finite measure and $\Om$ is an optimal set for \eqref{sop}. Then there exists $L>0$ such that for every open set $U$ containing $\Dr$ we have $\Om\subset U+B_L(0)$. In particular if $\Dr$ is bounded, then $\Om$ is bounded.
\end{prop}

\begin{proof}
It is enough to consider only the case of a quasi-connected set $\Dr$ and to work with \eqref{op}.

Assume by contradiction that such $L$ does not exist. Then, there is a sequence $(x_n)_{n\ge1}\subset\Omega$ such that $d(x_n, U)\rightarrow +\infty$ and $|x_n-x_m|\ge2r_0$, when $n\neq m$.
Since $\Omega=\{u>0\}$, we have $u(x_n)>0,\forall n\ge1$ and so, by Lemma \ref{sopra}, there are constants $C_0>0$ and $0<r_0$ such that for each $r<r_0$, we have the bound
$$\|u\|_{L^{\infty}(B_r(x_n))}\ge C_0r.$$
For each $n\ge1$, consider $y_n\in B_r(x_n)$ such that 
$$u(y_n)\ge\frac{1}{2}C_0r.$$
Consider the function $\varphi_r(\cdot-y_n)$, as defined in \eqref{solwith1}. Then $u-\lambda\|u\|_{\infty}\varphi_r(\cdot-y_n)$ is subharmonic, since
$$\Delta(u-\lambda\|u\|_{\infty}\varphi_r(\cdot-y_n))\ge-\lambda u+\lambda\|u\|_{\infty}\ge0.$$
So, we have the inequalities
$$\begin{array}{ll}
&\ds\int_{B_r(y_n)}\big(u(x)-\lambda\|u\|_{\infty}\varphi_r(x-y_n)\big)\,dx\ge|B_r|\big(u(y_n)-\lambda\|u\|_{\infty}\varphi_r(x-y_n)\big)\\
&\ds\hskip6.3truecm\ge|B_r|\Big(\frac{C_0}{2}r-\lambda\|u\|_{\infty}r^2\varphi_1(0)\Big);\\
&\ds\int_{B_r(y_n)}u(x)\,dx\ge r^{2+d}\lambda\|u\|_{\infty}\|\varphi\|_{L^1}+|B_r|\Big(\frac{C_0}{2}r-\lambda\|u\|_{\infty}r^2\varphi_1(0)\Big).\\
\end{array}$$
Choose now $0<r<r_0$ small enough such that $\frac{C_0}{2}r-\lambda\|u\|_{\infty}r^2\varphi_1(0)>0$. Then there is a constant $c>0$, such that for all $n\ge1$
$$\int_{B_r(y_n)}u(x)\,dx\ge c.$$
The fact that the balls $B_r(y_n)$ are all disjoint contradicts the integrability of $u$.
\end{proof}

\begin{oss}\label{relim1}
The constant $c$, depends on $C_0, r_0$ and $|\Dr|$. In fact, the proof of the proposition above gives an estimate on the number of admissible points $x_n$. Therefore the value of $L$, can be estimated more explicitly. Since $\Om$ is open and connected, the value of $L$ depends only on $\lb_1(\Dr), |\Dr|$ and the dimension $d$ of the space.
\end{oss}

\subsection{Convex constraint does not imply convex optimal set}

In this section we will prove that the solution $\Omega$ of the optimization problem \eqref{lbk1} might not be convex even if the constraint $\Dr$ is convex. Consider the sequence of constraints $(D_n)_{n\ge1}$, where $D_n=(-\frac{1}{n},\frac{1}{n})\times(-1,1)$ and consider the sequence of bounded open sets $(\Omega_n)_{n\ge1}$ such that for each $n$ big enough, $\Omega_n$ is a solution of the shape optimization problem:
\be
\min\big\{\lambda_1(\Omega)\ :\ D_n\subset\Omega,\ \Omega\ \mbox{ quasi-open},\ |\Omega|=m\big\}.
\ee

\begin{prop}
For every $m<4/\pi$, there is $N>0$ such that $\Omega_n$ is not convex for all $n\ge N$.
\end{prop}

\begin{proof}
We begin with some observations on the optimal sets.
\begin{enumerate}
\item By a Steiner simmetrization argument, all the sets $\Omega_n$ are Steiner symmetric with respect to the axes $x$ and $y$ (in consequence, they are also star shaped sets).\\
\item For $n$ large enough, we consider the set $\Omega'_n=D_n\cup B^*(m-\frac{4}{n})$, where for any $a>0$, $B^*(a)$ denotes the ball centered in $0$ of measure $a$. By the optimality of $\Omega_n$, we have
$$\lambda_1(\Omega_n)\le\lambda_1(\Omega'_n)\le\lambda_1(B^*(m-\frac{4}{n})).$$
\end{enumerate}
By Theorem \ref{thcc}, there is a $\gamma$-converging subsequence still denoted by $(\Omega_n)_{n\ge1}$. Let $\Omega$ be the $\gamma$-limit of this subsequence. Then
$$\begin{array}{ll}
&\bullet\quad\ds\lambda_1(\Omega)\le\liminf_{n\to\infty}\lambda_1(\Omega_n)\le\liminf_{n\to\infty}\lambda_1(B^*(m-\frac{4}{n}))=\lambda_1(B^*(m));\\
& \\
&\bullet\quad\ds|\Omega|\le\liminf_{n\to\infty}|\Omega_n|=m.
\end{array}$$
Using the fact that the ball is the unique minimizer of $\lambda_1$ under a measure constraint, we obtain $\Omega=B^*(m)$.
Consider a ball $B'$ of center $(0,\sqrt{\frac{m}{\pi}}-\eps)$ and radius $\eps$ and a ball $B''$ of center $(0,-\sqrt{\frac{m}{\pi}}+\eps)$ and radius $\eps$. Then
$$\Omega_n\cap B'\xrightarrow[n\to\infty]{\gamma}\Omega\cap B'=B',\qquad
\Omega_n\cap B''\xrightarrow[n\to\infty]{\gamma}\Omega\cap B''=B''.$$
Then there is some $n$ large enough such that both sets $B'\cap\Omega_n$ and $B''\cap\Omega_n$ are non-empty, and $\Omega_n$ cannot be convex. In fact, if by contradiction $\Omega_n$ was convex, then we should have that the rombus $R$ with vertices $(-1,0),(0,(0,-\sqrt{\frac{m}{\pi}}+\eps),(1,0),(0,\sqrt{\frac{m}{\pi}}-\eps))$ is contained in $\Omega_n$. But
$$|R|=2(\sqrt{\frac{m}{\pi}}-\eps)> m$$
for $\eps$ small enough and $m\le4/\pi$, and this is a contradiction.
\end{proof}

\subsection{Lack of monotonicity}
We show here that in problem \eqref{lbk1} the optimal solutions are not monotone with respect to $m$, i.e. $m_1<m_2$ does not imply in general that $\Om_1 \subset \Om_2$ where $\Om_i$ is optimal with the constraint $m_i$. Similarly, in the penalized problem \eqref{sop}, the same lack of monotonicty occurs with respect to $\Lb$, i.e. $\Lb_1< \Lb_2$ does not imply in general that $\Om^1 \supset \Om^2$ where $\Om^i$ is optimal with the penalization $\Lb_i$. Here we consider only the case of penalization, since the first one follows as a consequence, taking $m_1=|\Om^2|$ and $m_2=|\Om^1|$.

Let us consider in $\R^2$ the internal constraint $\Dr$ of the form $\Dr=B_{1/2}(0)\cup R^{\eps, \eta}$ where $R^{\eps, \eta}$ is the rectangle
$({\eta},0)+(-\frac{\eps}{2},\frac{\eps}{2})\times(-\frac{1}{2\eps},\frac{1}{2\eps})$. The parameters $\eps, \eta$ will be fixed later.

Note that $\frac{\pi}{4}=|B_{1/2}(0)|<|R^{\eps, \eta}|=1$ and that $\lb_1(B_{1/2}(0))< \lb_1(R^{\eps, \eta})$ for $\eps$ small enough. As well, we notice that the distance between $B_{1/2}(0)$ and $R^{\eps, \eta}$ tends to $+\infty$ as $\eta\to+\infty$. Following Remark \ref{relim1} for every $\Lb$ and $\eps>0$, there exists $\eta$ large enough such that every solution $\Om$ of \eqref{sop} satisfies one of the following two possibilities:
\begin{itemize}
\item[(A)] $\Om= B\cup R^{\eps, \eta}$, where $B$ is a ball containing $B_{1/2}(0)$ and disjoint from $R^{\eps, \eta}$;
 \item[(B)] $\Om=B_{1/2}(0) \cup A$, where $A$ is a connected open set containing $R^{\eps, \eta}$ and disjoint from $B_{1/2}(0)$.
\end{itemize}

\begin{lemma}\label{born005}
Let $\Lb>0$ be fixed, let $\Om_\eps$ be a solution of the problem 
$$\min\big\{\lb_1(\Om)+\Lb|\Om|\ :\ \Om\supset R^{\eps,0}\big\},$$
and let $B$ be a ball solving
$$\min\big\{\lb_1(\Om)+\Lb|\Om|\ :\ \Om\subset\R^2\big\}.$$
Then we have
$$\lb_1(B)=\lim_{\eps\to0}\lb_1(\Om_\eps),\qquad|B|+1=\lim_{\eps\to0}|\Om_\eps|.$$
\end{lemma}

\begin{proof}
By Steiner symmetrization along both axes, the sets $\Om_\eps$ are Steiner symmetric, and so star shaped. Therefore the sets $\Om_\eps$ fulfill a uniform exterior segment condition which, together with the compactness Theorem \ref{comgam}, is enough (see \cite[Chapter 4]{bubu05}) to give that $\Om_\eps$ $\g_{loc}$ converges to some open set $\Om$.

We first notice that 
\begin{equation}\label{born00}
\lb_1(\Om_\eps)+ \Lb |\Om_\eps|\le \lb_1(B_1(0))+\Lb |B_1(0)| + \Lb:=c,
\end{equation}
which gives that both measure of $\Om_\eps$ and $\lb_1(\Om_\eps)$ are uniformly bounded. Because of that and of the Steiner symmetrization above, all $\Om_\eps$ are contained in the set 
\begin{equation}\label{born001}
\{(x,y)\ :\ |xy|\le c\}.
\end{equation}

From the properties of the $\g_{loc}$-convergence, for every ball $B_R(0)$ we have that
$$|\Om\cap B_R(0)|\le\liminf_{\eps\to0}|\Om_\eps\cap B_R(0)|.$$
Since
$$\liminf_{\eps\to0}|\Om_\eps\cap B_R(0)|\le\liminf_{\eps\to0}|\Om_\eps|-1$$
we get
\begin{equation}\label{born004}
|\Om|+1\le\liminf_{\eps\to0}|\Om_\eps|.
\end{equation}

We prove now that 
\begin{equation}\label{born002}
\lb_1(\Om)\le\liminf_{\eps\to0}\lb_1(\Om_\eps).
\end{equation}
Let $u_\eps$ be the first normalized eigenfunction on $\Om_\eps$. By the concentration compactness principle, we may have: compactness, vanishing or dichotomy. The vanishing is ruled out by the fact that in this case we would have $\lb_1(\Om_\eps) \rightarrow +\infty$, which contradicts \eqref{born00}. The dichotomy is ruled out too, by the following argument. Let $u_\eps^i$, $i=1,2$ be the two sequences provided by the dichotomy. From the concentration compactness principle, at least one sequence of quasi-open sets $\{u_\eps^i >0\}$ has a distance from the origin going to $+\infty$. In the same time $\lb_1(\{u_\eps^i >0\}) $ are equibounded. This is in contradiction with the inclusion \eqref{born001}. Therefore the compactness occurs, i.e. $u_\eps(\cdot + y_\eps) $ converges strongly in $L^2(\R^2)$ to some function $u\in H^1_0(\Om)$. Again, by Steiner symmetrization the vectors $y_\eps$ can be taken equal to $0$. Consequently \eqref{born002} is achieved.

Taking test domains of the form $B\cup R^{\eps, 0}$ with $B\cap R^{\eps, 0} =\emptyset$ we have that
$$\lb_1(B)+\Lb|B|+\Lb\ge \lb_1(\Om_\eps)+\Lb|\Om_\eps| $$
and passing to the limit
$$\lb_1(B)+\Lb|B|+\Lb\ge \lb_1(\Om)+\Lb|\Om|+\Lb.$$
Using the optimality of the ball $B$ we get $\Om=B$ and inequalities \eqref{born004}-\eqref{born002} become equalities.
\end{proof}

Let us fix $\Lb_2$ such that a global solution of 
$$\min\big\{\lb_1(\Om)+\Lb_2|\Om|\ :\ \Om\subset\R^2\big\}$$
is the ball $B_1(0)$. Then for $\eps$ small enough given by Lemma \ref{born005} and for $\eta$ large enough given by Remark \ref{relim1} the solution of \eqref{sop} with $\Lb_2$ is
$$\Om_\eps^2=B_1(0)\cup R^{\eps,\eta}.$$
Indeed, from Lemma \ref{born005}, for $\eps$ small enough we have
$$\lb_1(\Om_\eps)+\Lb_2|\Om_\eps|+\Lb_2|B_{1/2}(0)|>\lb_1(B_1(0))+\Lb_2|B_1(0)|+\Lb_2|R^{\eps,\eta}|,$$
so situation (A) occurs. 

For the $\eps$ fixed above, take $\Lb_1$ small enough such that a ball $B'$ containing $R^{\eps,0}$ is a global minimizer for
$$\min\big\{\lb_1(\Om)+\Lb_1|\Om|\ :\ \Om\subset\R^2\big\}.$$
Then we are in situation (B) since $|B_{1/2}(0)|<|R^{\eps,0}|$. This concludes our argument since no monotonicity may occur.

\subsection{The optimal set has finite perimeter}

\begin{prop}
Assume that $\Dr$ is open and connected in \eqref{sop}. Then the perimeter in $\R^d\setminus\overline\Dr$ of an optimal set $\Om$ is finite.
\end{prop}

\begin{proof}
Let $u$ be a normalized eigenfunction associated to an optimal set $\Om$. Since $\Dr$ is connected we have $ \Omega=\{u>0\}$.
Consider the set $\Omega_{\eps}=\Dr\cup \{u>\eps\}$. By the optimality of $\Omega$ we have, using the fact that for $\eps$ small $|\{u\le\eps\}|\le|\{u>\eps\}|$,
\be
\begin{array}{ll}
\lambda_1(\Omega)+|\Omega|&\ds\le\lambda_1(\Omega_{\eps})+|\Omega_{\eps}|\le\frac{\int|\nabla (u-\eps)^{+}|^2\,dx}{\int|(u-\eps)^{+}|^2\,dx}+|\Omega_{\eps}|\\
\\
&\ds\le\frac{\int_{\{u>\eps\}}|\nabla u|^2\,dx}{\int_{\{u>\eps\}} (u^2-2\eps u+\eps^2)\,dx}+|\Omega_{\eps}|\le\frac{\int_{\{u>\eps\}}|\nabla u|^2\,dx}{1-2\eps\int_{\{u>\eps\}}u(x)\,dx}+|\Omega_{\eps}|\\
\\
&\ds\le\int_{\{u>\eps\}}|\nabla u|^2\,dx+\Big(\frac{1}{1-2\eps\int u(x)\,dx}-1\Big)\lambda_1(\Omega_{\eps})+|\Omega_{\eps}|\\
\\
&\ds\le\int_{\{u>\eps\}}|\nabla u|^2\,dx+2\eps\int u(x)\,dx\frac{1}{1-2\eps\int u(x)\,dx}\lambda_1(\Omega)+|\Omega_{\eps}|.
\end{array}
\ee
Then we have a constant $C$ depending on $u$ such that for $\eps$ small enough we have the inequalities
\be
\begin{array}{ll}
\eps C&\ds\ge\int_{\{0<u\le\eps\}}|\nabla u|^2\,dx+|\Omega\setminus\Omega_{\eps}|\\
\\
&\ds\ge\int_{\{0<u\le\eps\}\setminus\Dr}|\nabla u|^2\,dx+|\Omega\setminus\Omega_{\eps}|\\
\\
&\ds\ge\int_{\{0<u\le\eps\}\setminus\overline {\Dr}}|\nabla u|^2\,dx+|\{0<u\le\eps\}\setminus\overline {\Dr}|\\
\\
&\ds\ge\frac{1}{|\{0<u\le\eps\}\setminus\overline{\Dr}|}\Big(\int_{\{0<u\le\eps\}\setminus\overline{\Dr}}|\nabla u|\,dx\Big)^2+|\{0<u\le\eps\}\setminus\overline{\Dr}|.
\end{array}
\ee
Thus, there exists a constant $C$ independent on $\eps$ such that
\be\label{per}
\int_{\{0<u\le\eps\}\setminus\overline{\Dr}}|\nabla u(x)|\,dx\le C\eps.
\ee
By the co-area formula

\be\label{coarea}
\frac{1}{\eps}\int_0^{\eps}P\left(\{u>t\};\R^d\setminus\overline\Dr\right)\,dt\le C,
\ee
for each $\eps>0$ small enough. Then, there is a sequence $(\eps_n)_{n\ge1}$ converging to $0$ such that 
$$P\left(\{u>\eps_n\};\R^d\setminus\overline\Dr\right)\le C.$$
Passing to the limit we have 
$$P\left(\{u>0\};\R^d\setminus\overline\Dr\right)\le C.$$
\end{proof}

\begin{oss}
The regularity of the free parts of the boundary is the same as in \cite[Theorem 1.2]{brla}, being independent on the fact that the inclusion constraint is internal or external.
\end{oss}

\section{Open problems and complements}\label{scomm}

We give a list of some open problems that arose during the work on this article. We denote by $\Omega(\Dr,m)$ a quasi-open set of Lebesgue measure $m$, which solves the shape optimization problem \eqref{lbk1}.
\begin{enumerate}
\item Is there some $\eps>0$, such that for every $0<m<\eps$, the set $\Omega(\Dr,|\Dr|+m)$ is unique? Note that this is certainly not true when $m$ is large, since for a bounded $\Dr$ any ball of measure $m$ and containing $\Dr$ is a solution.

\item If $m'>m$, is there an optimal set $\Omega(\Dr,m')$ containing the optimal set $\Omega(\Dr,m)$? Note that the symmetric statement (if $m'<m$, then for each optimal set $\Omega(\Dr,m)$, there is an optimal set $\Omega(\Dr,m')\subset \Omega(\Dr,m)$) is false. Indeed, take for instance $\Dr$ the unit square in $\R^2$ centered in $0$ and $m'=\frac{\pi}{2}<m$. Then $\Omega(\Dr,m')$ is the ball centered at $0$ with radius $\frac{1}{\sqrt{2}}$ and $\Omega(\Dr,m)$ is any ball of radius $\sqrt{\frac{m}{\pi}}$. Clearly, there are balls $\Omega(\Dr,m)$ which do not contain $\Omega(\Dr,m')$.

\item Let $\Dr$ be an open convex set such that for every $m\ge|\Dr|$ there exists a convex solution to the shape optimization problem \eqref{lbk1}. Is it true that then $\Dr$ is a ball?

\item A interesting problem, similar to \eqref{lbk1}, is given by the minimization of the energy integral functional
$$E(\O)=\int-w_\O(x) dx.$$
We can repeat in this case all the arguments above, obtaining similar existence, boundedness and regularity results. In particular, working with the energy functional simplifies the analysis of Proposition \ref{open}, obtaining that optimal sets are open, even without the quasi-connectedness assumption on $\Dr$.

\item Is it true that for every quasi-open set $\Dr$ with finite measure, problem \eqref{lbkk} has a solution? \end{enumerate}

\begin{ack}
The work of Dorin Bucur is part of the project ANR-09-BLAN-0037 {\it Geometric analysis of optimal shapes (GAOS)} financed by the French Agence Nationale de la Recherche (ANR). The work of Giuseppe Buttazzo and Bozhidar Velichkov is part of the project 2008K7Z249 {\it Trasporto ottimo di massa, disuguaglianze geometriche e funzionali e applicazioni} financed by the Italian Ministry of Research.
\end{ack}


\bigskip
{\small\noindent
Dorin Bucur:
Laboratoire de Math\'ematiques (LAMA),
Universit\'e de Savoie\\
Campus Scientifique,
73376 Le-Bourget-Du-Lac - FRANCE\\
{\tt dorin.bucur@univ-savoie.fr}\\
{\tt http://www.lama.univ-savoie.fr/$\sim$bucur/}

\bigskip\noindent
Giuseppe Buttazzo:
Dipartimento di Matematica,
Universit\`a di Pisa\\
Largo B. Pontecorvo 5,
56127 Pisa - ITALY\\
{\tt buttazzo@dm.unipi.it}\\
{\tt http://www.dm.unipi.it/pages/buttazzo/}

\bigskip\noindent
Bozhidar Velichkov:
Scuola Normale Superiore di Pisa\\
Piazza dei Cavalieri 7, 56126 Pisa - ITALY\\
{\tt b.velichkov@sns.it}

\end{document}